\documentclass{ifacconf}
\usepackage{graphicx}      
\usepackage{natbib}        
\usepackage{amsmath,amssymb}
\usepackage{amsfonts}
\usepackage{mathtools}
\usepackage[new]{old-arrows}
\usepackage{array}
\usepackage{enumerate}
\usepackage[usenames,dvipsnames]{color}
\usepackage{colortbl, xcolor}
\usepackage{subfigure}
\usepackage{tikz}
\usetikzlibrary{positioning,arrows}
\allowdisplaybreaks

\newcommand{\R}{\mathbb{R}}
\newcommand{\cM}{\mathcal M}

\begin{document}
	\begin{frontmatter}
		
		\title{Approximation of a compound-exchanging cell by a Dirac point} 
		
		\author[First]{Xiao Yang} 
		\author[Second]{Qiyao Peng} 
		\author[First]{Sander C. Hille}
		
		\address[First]{Mathematical Institute, Leiden University, 
			Leiden, the Netherlands \emph{(}e-mail: x.yang@math.leidenuniv.nl; shille@math.leidenuniv.nl\emph{)}}
		\address[Second]{Bioinformatics, Computer Science Department, Vrije Universiteit Amsterdam, Amsterdam, the Netherlands \emph{(}e-mail: q.peng@vu.nl\emph{)}}
		
		\begin{abstract}
			Communication between single cells or higher organisms by means of diffusive compounds is an important phenomenon in biological systems. Modelling therefore often occurs, most straightforwardly by a diffusion equation with suitable flux boundary conditions at the cell boundaries. Such a model will become computationally inefficient and analytically complex when there are many cells, even more so when they are moving. We propose to consider instead a point source model. Each cell is virtually reduced to a point and appears in the diffusion equation for the compound on the full spatial domain as a singular reaction term in the form of a Dirac delta `function' (measure) located at the cell's centre. In this model, it has an amplitude that is a non-local function of the concentration of compound on the (now virtual) cell boundary. We prove the well-posedness of this particular parabolic problem with non-local and singular reaction term in suitable Sobolev spaces. We show for a square bounded domain and for the plane that the solution cannot be $H^1$-smooth at the Dirac point. Further, we show a preliminary numerical comparison between the solutions to the two models that suggests that the two models are highly comparable to each other.
		\end{abstract}
		
		\begin{keyword}
			Cell communication, Diffusion equation, Numerical analysis, Heat kernel, Dirac measure, Well-posedness 
		\end{keyword}
		
	\end{frontmatter}
	
	\section{Introduction}
	Cells communicate with each other to coordinate and accomplish tasks that they cannot do alone, for instance, the development of tissues, organs and whole organisms, and regulating the cell behaviours \citep{Alberts2002molecular, Neitzel2014Cell}. One means of cell-to-cell communication is via diffusive compounds (e.g. signalling molecules) released by the cell into the extracellular matrix (ECM), then binding to receptors of other cells or being taken up \citep{Muller_2011}. One may think of the secretion and detection of cytokines in immune response, the spreading of pheromones or signalling by means of ethylene among plants.
	In other words, cells (or higher organisms) do not only secret but may also `consume' diffusive compounds.

	Mathematically, we can model the concentration of the compounds over the domain by the diffusion equation. Depending on the scale and for the sake of computational efficiency, one can decide whether to separate the intra- and extracellular environment. \citet{HMEvers2015, Peng2023, Peng2024-fv} compared the solutions to the two models from analytical and numerical perspectives with different biological assumptions. In this paper, the names of the two models studied were inherited from \cite{Peng2023,Peng2024-fv}, namely, the \textit{spatial exclusion model} and the \textit{point source model}. Here, however, we consider both secretion and uptake: the secretion occurs with a predefined flux density and the uptake is proportional to the local concentration with a non-negative constant rate. Hence, the flux density in the spatial exclusion model and the intensity of the Dirac point in the point source model are modified accordingly from \cite{Peng2023}.
	
\begin{figure}[h!] \label{fig:schematic computational domain}
		\centering
		
		
		\tikzset {_vb1flt2vq/.code = {\pgfsetadditionalshadetransform{ \pgftransformshift{\pgfpoint{45 bp } { -74.5 bp}}\pgftransformscale{1 }  }}}
		\pgfdeclareradialshading{_p0z4j2a23}{\pgfpoint{-56bp}{72bp}}{rgb(0bp)=(0.84,1,0.76);
			rgb(8.482142857142858bp)=(0.84,1,0.76);
			rgb(25bp)=(1,1,0.91);
			rgb(400bp)=(1,1,0.91)}
		
		
		\tikzset {_iwr48r03h/.code = {\pgfsetadditionalshadetransform{ \pgftransformshift{\pgfpoint{45 bp } { -74.5 bp }  }  \pgftransformscale{1 }  }}}
		\pgfdeclareradialshading{_fqmewtzsb}{\pgfpoint{-56bp}{72bp}}{rgb(0bp)=(0.71,1,0.91);
			rgb(8.482142857142858bp)=(0.71,1,0.91);
			rgb(25bp)=(1,1,0.91);
			rgb(400bp)=(1,1,0.91)}
		
		
		\tikzset {_fvgw89tn6/.code = {\pgfsetadditionalshadetransform{ \pgftransformshift{\pgfpoint{-99 bp}{ 99 bp}}  \pgftransformscale{1.32}}}}
		\pgfdeclareradialshading{_bae41fu94}{\pgfpoint{80bp}{-80bp}}{rgb(0bp)=(1,1,1);
			rgb(0bp)=(1,1,1);
			rgb(25bp)=(1,1,1);
			rgb(400bp)=(1,1,1)}
		\tikzset{every picture/.style={line width=0.75pt}} 
		\begin{tikzpicture}[x=0.75pt,y=0.75pt,yscale=-0.3,xscale=0.3,every node/.style={scale=0.42}]
		
		\path  [shading=_fqmewtzsb,_iwr48r03h] (500,56) -- (854,56) -- (854,393) -- (500,393) -- cycle ; 
		\draw   (500,56) -- (854,56) -- (854,393) -- (500,393) -- cycle ; 
		
		\draw  [dash pattern={on 2.5pt off 2.5pt}] (543,254.5) .. controls (543,215.56) and (574.56,184) .. (613.5,184) .. controls (652.44,184) and (684,215.56) .. (684,254.5) .. controls (684,293.44) and (652.44,325) .. (613.5,325) .. controls (574.56,325) and (543,293.44) .. (543,254.5) -- cycle ;
		\path  [shading=_fqmewtzsb,_iwr48r03h] (29,56) -- (383,56) -- (383,393) -- (29,393) -- cycle ; 
		\draw   (29,56) -- (383,56) -- (383,393) -- (29,393) -- cycle ; 
		
		\path  [shading=_bae41fu94,_fvgw89tn6] (72,254.5) .. controls (72,215.56) and (103.56,184) .. (142.5,184) .. controls (181.44,184) and (213,215.56) .. (213,254.5) .. controls (213,293.44) and (181.44,325) .. (142.5,325) .. controls (103.56,325) and (72,293.44) .. (72,254.5) -- cycle ; 
		\draw   (72,254.5) .. controls (72,215.56) and (103.56,184) .. (142.5,184) .. controls (181.44,184) and (213,215.56) .. (213,254.5) .. controls (213,293.44) and (181.44,325) .. (142.5,325) .. controls (103.56,325) and (72,293.44) .. (72,254.5) -- cycle ; 
		
		\draw    (97,201) -- (132.75,245.44) ;
		\draw [shift={(134,247)}, rotate = 231.19] [color={rgb, 255:red, 0; green, 0; blue, 0 }  ][line width=0.75]    (10.93,-3.29) .. controls (6.95,-1.4) and (3.31,-0.3) .. (0,0) .. controls (3.31,0.3) and (6.95,1.4) .. (10.93,3.29)   ;
		\draw  [fill={rgb, 255:red, 0; green, 0; blue, 0 }  ,fill opacity=1 ] (611,255.5) .. controls (611,254.67) and (611.67,254) .. (612.5,254) .. controls (613.33,254) and (614,254.67) .. (614,255.5) .. controls (614,256.33) and (613.33,257) .. (612.5,257) .. controls (611.67,257) and (611,256.33) .. (611,255.5) -- cycle ;
		
		\draw (386,84.4) node [anchor=north west][inner sep=0.75pt]    {$\partial \Omega$};
		\draw (184,180.4) node [anchor=north west][inner sep=0.75pt]    {$\partial\Omega_{C}$};
		\draw (131,277.4) node [anchor=north west][inner sep=0.75pt]    {$\Omega_C$};
		\draw (137,229.4) node [anchor=north west][inner sep=0.75pt]    {$\mathbf{n}$};
		\draw (4,227.4) node [anchor=north west][inner sep=0.75pt]  [font=\normalsize]  {$\Omega$};
		\draw (265,343.4) node [anchor=north west][inner sep=0.75pt]    {$\tilde\Omega:=\Omega\setminus\Omega_C$};
		\draw (600,263.4) node [anchor=north west][inner sep=0.75pt]  [font=\normalsize]  {$\delta _{\boldsymbol{x}_C}$};
		\draw (475,227.4) node [anchor=north west][inner sep=0.75pt]  [font=\normalsize]  {$\Omega$};
		\draw (857,84.4) node [anchor=north west][inner sep=0.75pt]  {$\partial\Omega$};
		\draw (655,180.4) node [anchor=north west][inner sep=0.75pt]    {$\partial\Omega_C$};
		\end{tikzpicture}
		\caption{\it A schematic presentation of the computational domain  of the spatial exclusion model \emph{(}left\emph{)} and the point source model \emph{(}right\emph{)}. An example of diffusive spread of compound is indicated \emph{(}blue\emph{)}.} 
	\end{figure}
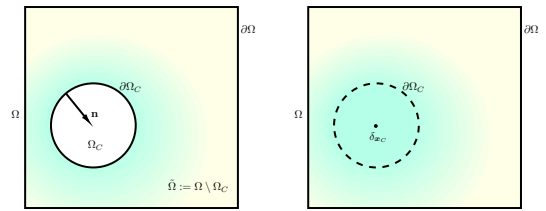

	In the \textit{spatial exclusion model}, a cell is modelled as a circle $\Omega_C$ of radius $R$ centered at $\boldsymbol{x}_C$ and is strictly embedded in a square domain $\Omega$ (see Fig. \ref{fig:schematic computational domain}, left). Diffusive compounds are only present in the extracellular environment or ECM, i.e.  $\tilde{\Omega}:=\Omega\setminus\Omega_C$. They are released from the cell boundary, with prescribed flux density $\phi(\boldsymbol{x},t)$. That is, for the computational domain, the intracellular environment is excluded from the square domain $\Omega$. The sign of $\phi(\boldsymbol{x},t)$ determines whether the cell secrets (positive) or takes up (negative) the compounds. 
	
	This process is modelled by the following diffusion partial differential equation together with mixed boundary and initial conditions:
	\begin{equation*}\label{DS}
	(\rm BVP_S)\quad\begin{cases}
	\frac{\partial u_S}{\partial t}=D\Delta u_S, & \text{in }\tilde\Omega, \\
	D\nabla u_S\cdot\mathbf{n} = \phi(\boldsymbol{x},t) - au_S, & \text{on }\partial\Omega_C, \\
	D\nabla u_S\cdot\mathbf{n} = 0, & \text{on } \partial\Omega, \\
	u_S(\boldsymbol{x},0) = u_0(\boldsymbol{x}), & \text{in $\tilde\Omega$}, \\
	\end{cases}
	\end{equation*}
	where $\phi(\boldsymbol{x},t)\geqslant 0$ and $a\geqslant0$ is constant. It indicates that the cell secrets the compounds over the cell boundary and takes up the compounds only when there is any present. The spatial exclusion model provides an accurate reflection of physical reality, but can become computationally `heavy' when many cells are involved. Moreover, in scenarios with moving cells -- which we need to consider in future work -- implementation of the time-varying computational domain will be difficult and require e.g. time-consuming remeshing in time-steps in a finite-element methods (FEM) implementation, rendering it (too) slow.
	
	Hence, the \textit{point source model} is favoured for more complicated scenarios, since the computational domain does not change. There, the compound distribution $u_P$ is defined in the entire spatial domain $\Omega$, substituting the cell by a Dirac measure $\delta_{\boldsymbol{x}_C}$ concentrated at the centre of the cell $\Omega_C$ (see Fig. \ref{fig:schematic computational domain}, right) with amplitude 
	\begin{equation}\label{eq:def Psi}
	\Psi[u_P](t):=\oint_{\partial\Omega_C}\phi(\boldsymbol{x},t) - au_P(\boldsymbol{x},t)d\Gamma(\boldsymbol{x}).
	\end{equation} 
	Here, $\Gamma$ is the measure on $\partial\Omega_C$ such that the Divergence Theorem holds on $\Omega_C$ without a constant.
	The point source model in the full open bounded domain $\Omega$ is then given by
	\begin{equation*}
	(\rm BVP_P)\quad\begin{cases}
	\frac{\partial u_P}{\partial t}= D\Delta u_P +\Psi[u_P](t)\,\delta_{\boldsymbol{x}_C},& \text{in }\Omega, \\
	D\nabla u_P\cdot\mathbf{n} = 0, & \text{on } \partial\Omega, \\
	u_P(\boldsymbol{x},0) = \tilde{u}_0(\boldsymbol{x}), & \text{in $\Omega$}.
	\end{cases}
	\end{equation*}
	The amplitude $\Psi$ for the Dirac measure has been defined such that the total amount of compound that would have been exchanged over the (now virtual) cell boundary per unit time, is exchanged through the Dirac point $\boldsymbol{x}_c$ per unit time. This approach follows  \citet{Peng2023}. Note that $\Psi$ depends non-locally on the function $u_P$.

	We shall first address in Section \ref{sec:well_posedness} the question of how to interpret equation $(\rm BVP_P)$ that involves a Dirac measure. The measure prevents viewing it pointwise from the classical perspective. Next, we consider for each model whether the Cauchy problem in that specific notion of solution is well-posed in particular function space(s). Then in Section \ref{sec:numerical} we provide preliminary numerical comparisons of the solution of $(\rm BVP_S)$ and $(\rm BVP_P)$. Section \ref{sec:conclusion} delivers conclusions from the current work and an outlook on forthcoming research.

	\section{Well-posedness of both models}\label{sec:well_posedness}
	
	For any final time $T>0$, $(\rm BVP_S)$ and $(\rm BVP_P)$ may both be interpreted in a distributional sense on the space-time domain $\overline{\tilde\Omega}\times [0,T]$ and $\overline{\Omega}\times[0,T]$, respectively. However, the space of distributions is too large. Moreover, the modeling imposes an evolution point of view as natural, in which a solution $u$ is considered as a particular distribution, namely a function from the time interval $I=[0,T]$ to some Banach function space on the spatial domain that is continuous, or at least locally integrable. Taking $D=1$ in this section, this yields the distributional weak form for $(\rm BVP_S)$:
	\[
	({\rm WF_S})\int_{\tilde{\Omega}\times I}(-\partial_t\varphi-\Delta\varphi)u_S d\Omega dt=\int_{\partial\Omega_C\times I}\varphi(\phi-au_S)d\Gamma dt
	\]
	for any smooth $\varphi$ on $\overline{\tilde\Omega}\times I$. The (distributional) weak form of $(\rm BVP_P)$ becomes
	$$({\rm WF_P})\;\int_{\Omega\times I}(-\partial_t\varphi-\Delta \varphi)u_P d\Omega dt=\int_{\Omega\times I}\varphi\Psi d\delta_{\boldsymbol x_C}dt$$
	for any smooth $\varphi$ on $\overline{\Omega}\times I$. 
	
	We will take `solution' in these weak senses. Spatial and time regularity will be made precise below. The provided regularity results may not be optimal: We give conditions sufficient for our purposes.
	
	We shall use standard notation and norms for $L^p$-spaces and Sobolev-Slobodeckij spaces $W^{k,p}(\Omega)$ and $H^k(\Omega)=W^{k,2}(\Omega)$: See \cite{Adams:2003} for details. The Boundary Trace Embedding Theorem for Sobolev spaces, (see e.g. \cite{Adams:2003}, Theorem 5.36, p.164) yields the existence of a boundary trace operator, which we further restrict to the part $\partial\Omega_C$ of $\partial\Omega$. The composition is here denoted by $\mathrm{Tr}$. It maps $W^{1,p}(\Omega)$ continuously into $L^p(\partial\Omega)$. When clear from context, we may omit $\mathrm{Tr}$ in $\mathrm{Tr}(u)$.
	
	
	

	\subsection{Well-Posedness of the spatial exclusion model}
	
	For the spatial exclusion model, a well-posedness result in $H^1(\tilde{\Omega})$ can be straightforwardly obtained from the extensive literature on parabolic problems, in particular the heat equation, since no singularity is involved. We follow exposition and results from \citet{Salsa2015partial}.
	
	
	
	
	\vskip 2mm
	
	\begin{thm}
		\label{thm:bvps}
		For every $T>0$, $u_0\in L^2(\tilde\Omega)$ and 
		\begin{equation}\label{eq:cond phi}  
		\phi\in L^2((0,T), H^1(\partial\Omega_C))\cap H^1((0,T),L^2(\partial\Omega_C)),
		\end{equation}
		$(\rm BVP_S)$ has a unique weak solution 
		\begin{equation} \label{eq:smoothness u}
		u_S\in L^2((0,T),H^1(\tilde\Omega))\cap H^1((0,T), H^1(\tilde{\Omega})^*).
		\end{equation}
		Moreover, $u_S$ depends continuously on initial condition $u_0$.
	\end{thm}
	\begin{pf}
		This can proven by the theorems in Section 10.3 and Section 10.4 of \citet{Salsa2015partial}. `Weak solution' is now stronger, implying $(\rm WF_S)$: 
		\[
		u_S\in 
		\bigl\{u\in L^2((0,T), H^1(\tilde{\Omega})\colon \partial_tu \in L^2((0,T),H^1(\tilde{\Omega})^*)\bigr\}
		\]
		satisfies for all $v\in L^2(\tilde{\Omega})$ and a.e. $t\in (0,T)$: 
		\[
		\langle \partial_t u_S(t), v\rangle + B(u(t),v) =0, 
		\]
		where $B$ is the bilinear form
		\[
		B(u,v):=\int_{\tilde\Omega}\nabla u\cdot\nabla v\ d\Omega+ \oint_{\partial\Omega_C} auv\ d\Gamma.
		\]
		H\"olders Inequality (applied twice) yields
		\[
		|B(u,v)|\leqslant \Vert u\Vert_{H^1}\Vert v\Vert_{H^1} + a\Vert u\Vert_{L^2(\partial\Omega_C)}\Vert v\Vert_{L^2(\partial\Omega_C)}.
		\]
		Hence, $B(\cdot, \cdot)$ is continuous on $H^1(\tilde{\Omega})$, because $\mathrm{Tr}: H^1(\tilde{\Omega})\to L^2(\partial\Omega_C)$ is continuous. 
		As $a\geqslant0$,
		$$B(u,u)\geqslant\Vert\nabla u\Vert_{L^2}^2=\Vert u\Vert_{H^1}^2-\Vert u\Vert_{L^2}^2.$$
		That is, $B(\cdot, \cdot)$ is $H^1$-$L^2$ weakly coercive. \citet{Salsa2015partial}, Theorem 10.6, then gives the result for $\phi=0$, with solution $\hat{u}_S$, say. For $\phi\neq 0$ as in \eqref{eq:cond phi}, Theorem 5.1 of \citet{HMEvers2015} yields a unique `particular solution' $\Phi$, as smooth as $u_S$ in \eqref{eq:smoothness u}, that solves $(\rm BVP_S)$ with $a=0$ and $u_0=0$. Then $u_S:=\hat{u}_S+\Phi$ is the desired solution.
	\end{pf}
	
	\subsection{Well-posedness of the point source model}
	Essentially, the argumentation in this section specializes to our particular setting, where the very general results obtained in \cite{Amann2001, AMQ:03}. In the latter, semilinear parabolic problems of the form
	\[
	\partial_t u = Au + F(u),\quad u(0)=u^0\in X,
	\]
	are studied, where $X$ is a Banach subspace of a suitable Sobolev space on $\Omega$ and $(A,\mathcal{D}(A))$ is the generator of a strongly continuous analytic semigroup $(e^{At})_{t\geqslant 0}$ on $X$. It is obtained by interpolation/extension methods from the Laplacian (in terms of distributional derivatives on the domain $W^{2,p}_B(\Omega)$. Here, subscript $B$ indicates that it is the closure in $W^{2,p}(\Omega)$ of the compactly supported smooth functions on $\overline{\Omega}$ that satisfy the Neumann boundary condition (cf. \cite{Amann2001}, Lemma 8). $F$ is a suitably regular map between $L^q_{\mathrm{loc}}((0,T), X)$ and $L^r_{\mathrm{loc}}((0,T), Y)$. The Banach space $Y$ is a subspace of a Sobolev space obtained by extrapolation from $X$, equipped with corresponding extension of $(S_t)_{t\geqslant 0}$.
	
	The key point in this section is to show that regularity parameters for the Sobolev spaces can be chosen in such a way that the conditions of \cite{AMQ:03} Theorem 3.4, p.1061 are satisfied, resulting in a well-posedness result. We need to restrict to the reflexive case and $1<p<2$, following \citet{AMQ:03}.
	\vskip 2mm
	
	\begin{thm}\label{thrm:well-posed point source}
		For any $p\in(1,2)$, $(\rm BVP_P)$ has a unique global solution $u_P\in C^0(\R_+,W^{1,p}(\Omega))$, provided $\tilde{u}_0\in W^{1,p}(\Omega)$ and $\phi\in L_{\mathrm{loc}}^\infty(\R_+\times\partial\Omega_C)$. Moreover, $u_P$ depends continuously on the initial condition $\tilde{u}_0$.
	\end{thm}
	
	Before getting to the proof, let us make two observations. Firstly, $W^{s,p}_B=W^{s,p}$ whenever $-1+\frac{1}{p}<s<1+\frac{1}{p}$ (cf. \cite{Amann2001}, p. 106). Secondly, Lemma 4.1(iii) in \cite{AMQ:03} implies that for any $\sigma\in [0,\frac{2}{p})\setminus\{1+\frac{1}{p}\}$, the Banach space $\cM(\overline{\Omega})$ of finite signed Radon measures on $\overline{\Omega}$ (with total variation norm) can be continuously embedded into $W^{\sigma-2,p}_B$.
	
	
	
	
	
	\begin{pf} {\it \emph{(}Theorem \ref{thrm:well-posed point source}\emph{)}.} In the notation of Theorem 3.4, \cite{AMQ:03}, we take $s=1$, $\rho=0$ and select $\displaystyle\sigma\in (1,\frac{2}{p})\setminus\{1+\frac{1}{p}\}$. Pick any $r>2/(\sigma-s)$, finite. Since $\phi\in L^\infty_{\mathrm{loc}}(\R_+\times\partial\Omega_C)$, it is in $L^r([0,T]\times\partial\Omega_C)$. Note that $W^{1,p}_B=W^{1,p}$. We will show that the map 
		\begin{gather*}
		F:C^0([0,T],W^{1,p}(\Omega))\to L^r([0,T],\cM(\overline{\Omega})),\\
		\qquad u\ \ \mapsto\ \ \bigl(\;t\mapsto \Psi[u](t)\delta_{\boldsymbol x_C}\bigr).
		\end{gather*}
		in ($\rm BVP_P$) is a bounded map and locally Lipschitz continuous. Then, because of the continuous linear embedding of $\cM(\overline{\Omega})$ into $W^{\sigma-2,p}_B(\Omega)$, $F$ satisfies the conditions of Theorem 3.4 and the stated results follow from the theorem, for initial condition $u^0\in W^{\sigma-2/r, p}_B(\Omega)$. Since $\sigma-\frac{2}{r}<s$ by assumption, a Sobolev Embedding Theorem and the first observation above yield that the result holds for any $u^0\in W^{1,p}(\Omega)=W^{1,p}_B(\Omega)\hookrightarrow W^{\sigma-2/r,p}_B$. The embedding is continuous.
		
		That $F$ is a bounded map, follows from
		\begin{align}
		\Vert F(u)\Vert_{L^r}&=\left(\int_0^T\left|\oint_{\partial\Omega_C}\phi(\boldsymbol{x},t)-au(\boldsymbol{x},t)d\Gamma\right|^rd t\right)^{\frac{1}{r}}\nonumber\\
		&\leqslant \Vert\phi\Vert_{L^r}+AT^{\frac{1}{r}}\textstyle\sup_{t\in[0,T]}\Vert u(\cdot,t)\Vert_{W^{1,p}}\nonumber\\
		&\leqslant A'(1+\textstyle\sup_{t\in[0,T]}\Vert u(\cdot,t)\Vert_{W^{1,p}}).\label{eq:bounded F}
		\end{align}
		For Lipschitz continuity, let $u,v\in C^0([0,T],W^{1,p}(\Omega))$. A similar computation as above yields
		\begin{align*}
		\Vert F(u)-F(v)\Vert_{L^r}^r&=a^r \int_0^T\left|\oint_{\partial\Omega_C}\mathrm{Tr}(u)-\mathrm{Tr}(v)\ d\Gamma\right|^rd t\\
		&\leqslant a^rC\int_0^T\Vert u(\cdot, t)-v(\cdot, t)\Vert_{W^{1,p}}^rd t\\
		&\leqslant a^r C T\textstyle\sup_{t\in[0,T]} \Vert u(\cdot, t)-v(\cdot, t)\Vert_{W^{1,p}}^r.
		\end{align*}
		Estimate \eqref{eq:bounded F} implies that the solution is in fact global.
	\end{pf}

	\subsection{ Singularity of solution at the Dirac point}
	
	The condition $p<2$ in Theorem \ref{thrm:well-posed point source} cannot be improved. The spatial regularity of the solution $u_P$ cannot become $H^1(\Omega)$, as we shall see below. 
	
	However, in the field of FEM, solutions are assumed to be the in the Sobolev space $H^1$ or $H^1_0$ on the spatial domain `by default'. So, some care and modification of the standard approach of FEM is needed. For $p<2$ the solution can still be approximated by typical families of basis functions, e.g. tent functions. Relevant results can be found in \cite{Scott1973}. 
	
	Space limitations prohibit a detailed analysis of the above regularity claim for general bounded domains $\Omega$. Instead, we consider specifically the square domain considered here, for which explicit computations can be made. Then we provide -- as further motivation -- a discussion of the heat equation on $\R^2$ with a single point source at the origin and zero initial data:  
	\begin{equation}
	\label{eq:2.1}
	\begin{cases}
	\frac{\partial \hat{u}}{\partial t}=\Delta \hat{u} + \delta_0, & \text{ in }\R^2,\\
	\hat{u}(\boldsymbol{x},0) = 0, & \text{ in }\R^2,
	\end{cases}
	\end{equation}
	because it suggests an approach to the general case. The explicit expression for the heat kernel for $\R^2$ can be replaced by appropriate estimates for the heat kernel and its gradient for Neumann boundary conditions that are available in the literature for suitably smooth and curved bounded domains.

	The theory of analytic semigroups (see e.g. \cite{Taira2016}) provides a framework to describe the time-evolution of System \eqref{eq:2.1}. Let $(e^{\Delta t})_{t\geqslant 0}$ denote the heat semigroup, which is given by convolution with the heat kernel on $\R^2$. Then, System \eqref{eq:2.1} can be solved as if it were an ODE in the Banach space completion of $\cM(\overline{\Omega})$ equipped with the Dudley or dual-bounded-Lipschitz norm, by means of the variation of constants formula:
	\begin{align*}
	\hat{u}(\boldsymbol{x},t) &= \ \bigl(\int_0^te^{\Delta(t-s)}\delta_0d s\bigr) (\boldsymbol{x})\\
	&=\ \int_0^t\int_{\R^2}\frac{1}{4\pi(t-s)}e^{-\frac{|\boldsymbol x-\boldsymbol x_0|^2}{4(t-s)}}\delta_0(d\boldsymbol x_0) d s\\
	&=\ \int_0^t\frac{1}{4\pi(t-s)}e^{-\frac{|\boldsymbol x|^2}{4(t-s)}}d s\\
	&=\ \int_{t^{-1}}^{+\infty}\frac{1}{4\pi \sigma}e^{-\frac{|\boldsymbol x|^2\sigma}{4}}d \sigma.
	\end{align*}
	
	
	\begin{prop}
		The solution $\hat{u}(\boldsymbol{x},t)$ to System \eqref{eq:2.1} has a logarithmic singularity around the origin for any $t>0$. Moreover, $\hat{u}(\boldsymbol{x}, t)\in L^2(\R^2)\setminus H^1(\R^2)$ for all $t>0$.
	\end{prop}
	\begin{pf}
		The solution $\hat{u}(\boldsymbol{x},t)$ is essentially the \textit{exponential integral}
		$$E_1(x):=\int_1^{+\infty}\frac{e^{-xs}}{s}d s.$$
		In fact, $\displaystyle \hat{u}(\boldsymbol{x},t) = \frac{1}{4\pi}E_1\left(\frac{|\boldsymbol{x}|^2}{4t}\right)$, and $E_1$ admits the expansion \citep{Abramowitz1948handbook}
		\[
		E_1(x)=-\log(x)-\gamma+\mathcal{O}(x),\qquad \mbox{as}\ x\to 0,
		\]
		where $\gamma$ is the Euler-Mascheroni constant. To determine whether $\hat{u}(\boldsymbol{x},t)$ is in $L^2(\R^2)$ for $t>0$, we therefore compute (up to a constant) in polar coordinates 
		\begin{align*}
		&\left(\int_0^{2\pi}\int_0^{+\infty} (E_1(r^2))^2 rd rd\theta\right)^{\frac{1}{2}}\\
		=&\sqrt{2\pi}\left(\int_0^{+\infty}\left(\int_1^{+\infty}\frac{e^{-r^2s}}{s}d s\right)^2rd r\right)^{\frac{1}{2}}\\
		\leqslant& \sqrt{2\pi}\int_1^{+\infty}\left(\int_0^{+\infty}\left(r^{\frac{1}{2}}\frac{e^{-r^2s}}{s}\right)^2d r\right)^{\frac{1}{2}}d s\\
		=&\frac{\sqrt{2\pi}}{4}\int_1^{+\infty}s^{-3/2}d s<\infty.
		\end{align*}
		In the inequality above we applied the Minkowski inequality for integrals. This shows that $\hat{u}(\cdot,t)\in L^2(\R^2)$ for all $t>0$. To determine whether it is in $H^1(\R^2)$, we continue with $\Vert\nabla \hat{u}\Vert_{L^2(\R^2)}$. This is tantamount to evaluating the following integral
		\begin{align*}
		\int_0^{+\infty}&\left(\frac{\partial}{\partial r}\int_1^{+\infty}\frac{e^{-r^2s}}{s}d s\right)^2rd r \\  
		& = \int_0^{+\infty}\left(\int_1^{+\infty}2re^{-r^2s}d s\right)^2rd r
		=  \int_0^{+\infty}\frac{4e^{-2r^2}}{r}d r,
		\end{align*}
		which is infinite. Hence, $\hat{u}(\cdot, t)$ is not in $H^1(\R^2)$ for any $t>0$. 
	\end{pf}
	
	It has been shown in \citet{HMEvers2015} -- in more generality -- that $\hat{u}(\boldsymbol{x},t)\in H^1(\R^2\setminus U)$ for all $t>0$, where $U$ is any open neighbourhood of the origin. Because of the mentioned kernel estimates, we expect a similar result to hold on $\Omega$: $u_P(\boldsymbol{x},t)\in H^1(\Omega\setminus U)$ for all $t>0$ for open neighbourhoods $U$ of $\boldsymbol{x}_C$. Thus, we can compare the solutions to the two models under $L^2$-norm and $H^1$-norm on $\tilde\Omega$.
	
	For a square domain, say $\Omega=(0,\pi)\times(0,\pi)$, one has a series expansion for the Neumann heat kernel at time $t$ (with $D=1$, still), with $\boldsymbol{x}_0=(x_0,y_0)$ and $\boldsymbol{x}=(x,y)$:
	\begin{gather*}
	h_t(\boldsymbol x,\boldsymbol x_0) =\frac{1}{\pi^2}+\sum_{m=1}^{\infty}\frac{2}{\pi^2}e^{-m^2t}\cos(mx)\cos(mx_0)\\
	\qquad\qquad +\ \sum_{m=1}^\infty\frac{2}{\pi^2}e^{-m^2t}\cos(my)\cos(my_0)\\
	+\sum_{m,n=1}^\infty \frac{4}{\pi^2}e^{-(m^2+n^2)t}\cos(mx)\cos(ny)\cos(mx_0)\cos(ny_0).
	\end{gather*}
	
	It allows us to make similar computations as was done for $\R^2$ for $\Omega$. Take $\phi\equiv 1$, $a=0$, $u_0\equiv 0$ with $\boldsymbol{x}_C=(\xi,\eta)$ as illustrative example. Then
	\begin{gather*}
	u_P(\boldsymbol x,t)\ =\ \int_0^t\int_{\Omega}h_{t-s}(\boldsymbol x,\boldsymbol x_0)d\delta_{(\xi,\eta)}(\boldsymbol{x_0}) d s=\frac{t}{\pi^2}+\\
	\sum_{m=1}^{\infty}\frac{2-2e^{-m^2t}}{\pi^2m^2}(\cos(mx)\cos(m\xi)+\cos(my)\cos(m\eta))+\\
	\sum_{m,n=1}^\infty\frac{4-4e^{-(m^2+n^2)t}}{\pi^2(m^2+n^2)}\cos(mx)\cos(m\xi)\cos(ny)\cos(n\eta).
	\end{gather*}
	
	\begin{prop}
		Let $u_P$ be the solution to $\rm (BVP_P)$ with $\phi\equiv 1$, $a=0$, $u_0\equiv 0$, then for any $t>0$, $u_P(\boldsymbol{x},t)\in L^2(\Omega)\setminus H^1(\Omega)$.
	\end{prop}
	
	\begin{pf}
		There is a bijection between $L^2(\Omega)$ and $\ell^2$ via Fourier series. In order to show $u_P$ is square-integrable, it is equivalent to show the Fourier coefficients are square-summable. As the numerators of all coefficients are bounded, it is enough to see the denominator. We know that $$\sum_{m=1}^\infty\frac{1}{m^4}=\zeta(4)=\frac{\pi^4}{90}<\infty.$$
		It remains to show the double infinite series
		$$\sum_{m,n=1}^\infty\frac{1}{(m^2+n^2)^2}<\infty.$$
		By integral test, it is equivalent to show
		$$\int_1^{+\infty}\int_1^{+\infty}\frac{1}{(x^2+y^2)^2}d xd y<\infty.$$
		Transform to polar coordinates.
		$$\int_1^{+\infty}\int_1^{+\infty}\frac{1}{(x^2+y^2)^2}d xd y\leqslant\int_0^{\frac{\pi}{2}}\int_1^\infty \frac{r}{r^4}d rd\theta<\infty.$$
		Hence, $u_P\in L^2(\Omega)$.

		For $\nabla u_P$, it suffices to show that
		$$\sum_{m,n=1}^\infty\frac{m^2}{(m^2+n^2)^2}=\infty.$$
		Indeed,
		\begin{align*}
		\sum_{m,n=1}^\infty\frac{m^2}{(m^2+n^2)^2} & \geqslant \sum_{m,n=1}^\infty\frac{m^2}{(n+m)^4}=\sum_{l=2}^\infty\sum_{m=1}^{l-1}\frac{m^2}{l^4}\\
		& = \sum_{l=2}^\infty\frac{l(l-1)(2l-1)}{6l^4}\\
		& =\sum_{l=2}^\infty\frac{1}{3l}-\sum_{l=2}^\infty\frac{1}{2l^2}+\sum_{l=2}^\infty\frac{1}{6l^3}.
		\end{align*}
		The latter expression contains two absolutely convergent series and one divergent series. Hence, it is divergent.
	\end{pf}

	\section{Numerical comparison}\label{sec:numerical}
	In this section, we present the results of numerical simulation and comparison of the two models on compounds-exchanging of the cell, i.e. secretion and uptake. An analytic comparison was initiated in \citet{HMEvers2015} on $\R^2$, without considering uptake, i.e. $a=0$. The non-regularity of the source term, i.e. being measure-valued, complicates the analysis already in that case. Here, the non-local character of the reaction term $\Psi[u_P](t)\delta_{\boldsymbol{x}_c}$ when $a>0$, is added to the challenge. This aspect of the comparison will be taken up in forthcoming work, after numerical investigations have been performed of which preliminary results are shown below.

	The simulation is conducted by FEM via \texttt{Python-FEniCSx (v0.7.3)} \citep{AlnaesEtal2015} for spatial discretization and backward (implicit) Euler method for time discretization. For the sake of maintaining the continuity of the weak solution, we approximate the Dirac measure by a Gaussian distribution with a very small variance \citep{Peng_2022_MATCOM}. The parameter values can be found in Table \ref{tab:parameter_values}, if there is no further specification. In this paper, we have focused mainly on the theoretical perspective, hence, we neglect dimensions for the numerical simulations. Moreover, for now, we only consider constant flux density defined in $\rm (BVP_S)$, that is, $\phi(\boldsymbol{x}, t) = \phi >0$.
	\begin{table}[h!]\footnotesize
		\centering
		\caption{\it Parameter values used in Section \ref{sec:numerical}}
		\begin{tabular}{m{1.5cm}m{1cm}m{5cm}}
			\hline
			\textbf{Parameter}  & \textbf{Value} & \textbf{Description} \\
			\hline
			$D$ & $1.0$ & Diffusion coefficient \\
			$\phi$ & $1.0$ & Secretion rate of compounds from the cell \\
			$a$ & $1$ & Uptake rate of compounds of the cell \\
			$\boldsymbol{x}_C$ & $(5,5)$ & Cell centre location\\
			$R$ & $0.25$ & Cell radius\\
			$L$ & $10$ & Side length of square computational domain $\Omega$\\
			$h$ & $0.2495$ & Average mesh size\\
			$\Delta t$ & $0.04$ & Time step\\
			$\varepsilon$ & $0.02$ & Variance in Gaussian distribution replacing the Dirac measure\\
			\hline
		\end{tabular}
		\label{tab:parameter_values}
	\end{table}

	When $a=0$, i.e. cell is only releasing the compounds and no uptake, then the model is the same as in \citet{Peng2023}, in which the singularity and the continuous linear increase of the $L^2$-norm in $\tilde{\Omega}$ have been observed. Moreover, there is a systematic time delay between the point source model and the spatial exclusion model, depending on the value of the diffusion coefficient -- compounds need time to diffuse from the cell centre to the cell boundary in the point source model. Then within the predefined time domain ($(0, 40.0)$ in \citet{Peng2023}), the difference between the solutions to the two models will stay constant after some time, while the total mass of the compounds in the computational domain increases over time.
	
\begin{figure}[h!]
		\centering
		\includegraphics[width=0.45\textwidth]{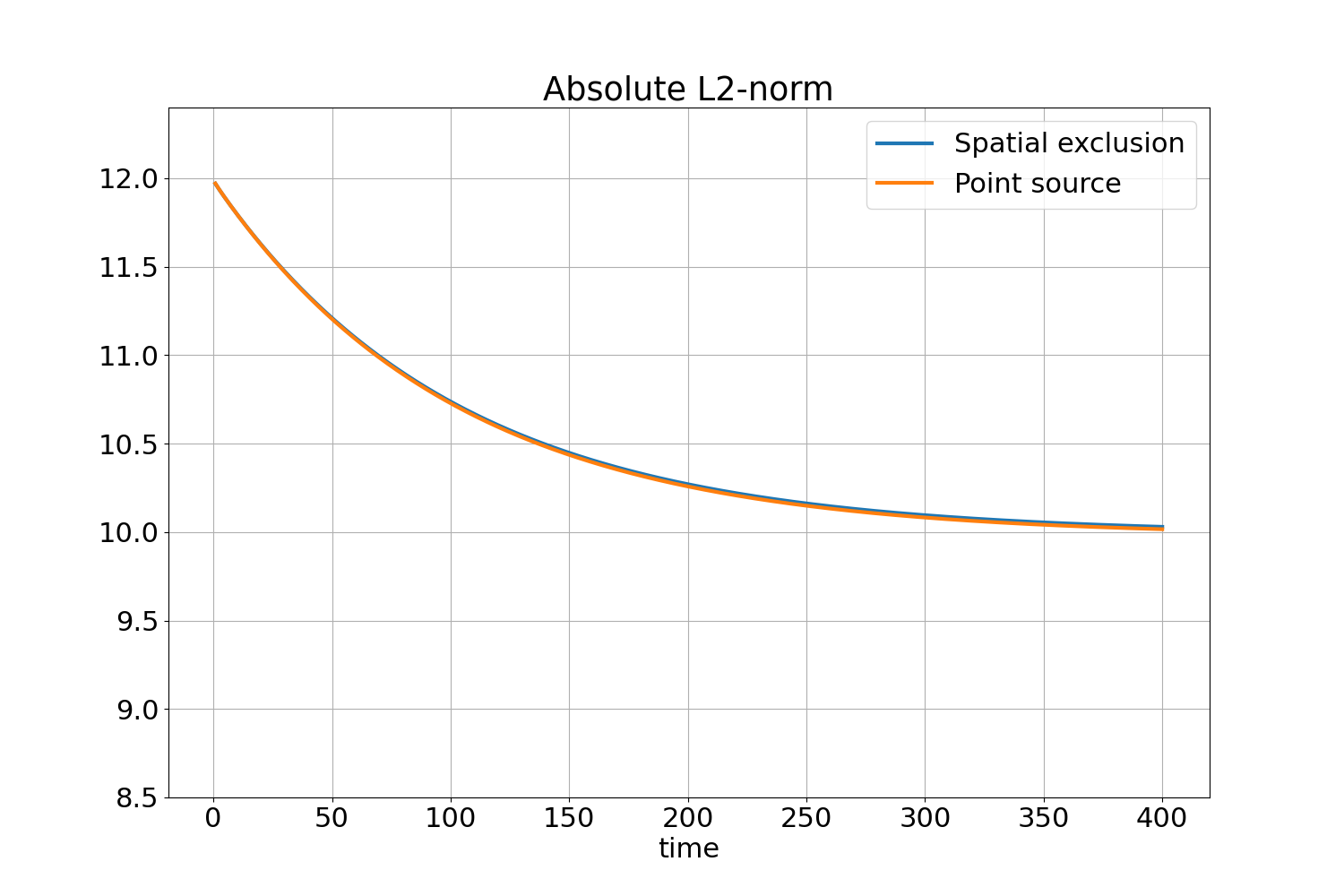}
		\caption{\it {Case when cell secrets and takes up the compounds \emph{(}i.e. both $\phi, a>0$\emph{)}: $L^2$-norm of solution to both models are shown, computed over the shared domain $\tilde{\Omega}$, with zero initial conditions \emph{(}i.e. $u_0 = \tilde{u}_0 = 0$\emph{)}. Both models will converge to a steady state. Blue and orange curves represent the spatial exclusion model and the point source model, respectively. Parameter values are taken from Table~\ref{tab:parameter_values}. }}
		\label{fig:steady}
	\end{figure}

	Considering uptake, i.e. $a>0$, the concentration of the compounds in both models, will reach a steady state, with a positive concentration of compounds in $\tilde{\Omega}$; see Figure \ref{fig:steady}. Thus, after some time, there is no appreciable exchange of compounds between the cell and the ECM.
	
	
	
	
	
	\begin{figure}[h!]
		\centering
		\subfigure[Only secretion: $\phi = 1.0, \ a = 0$]{
			\includegraphics[width=0.4\textwidth]{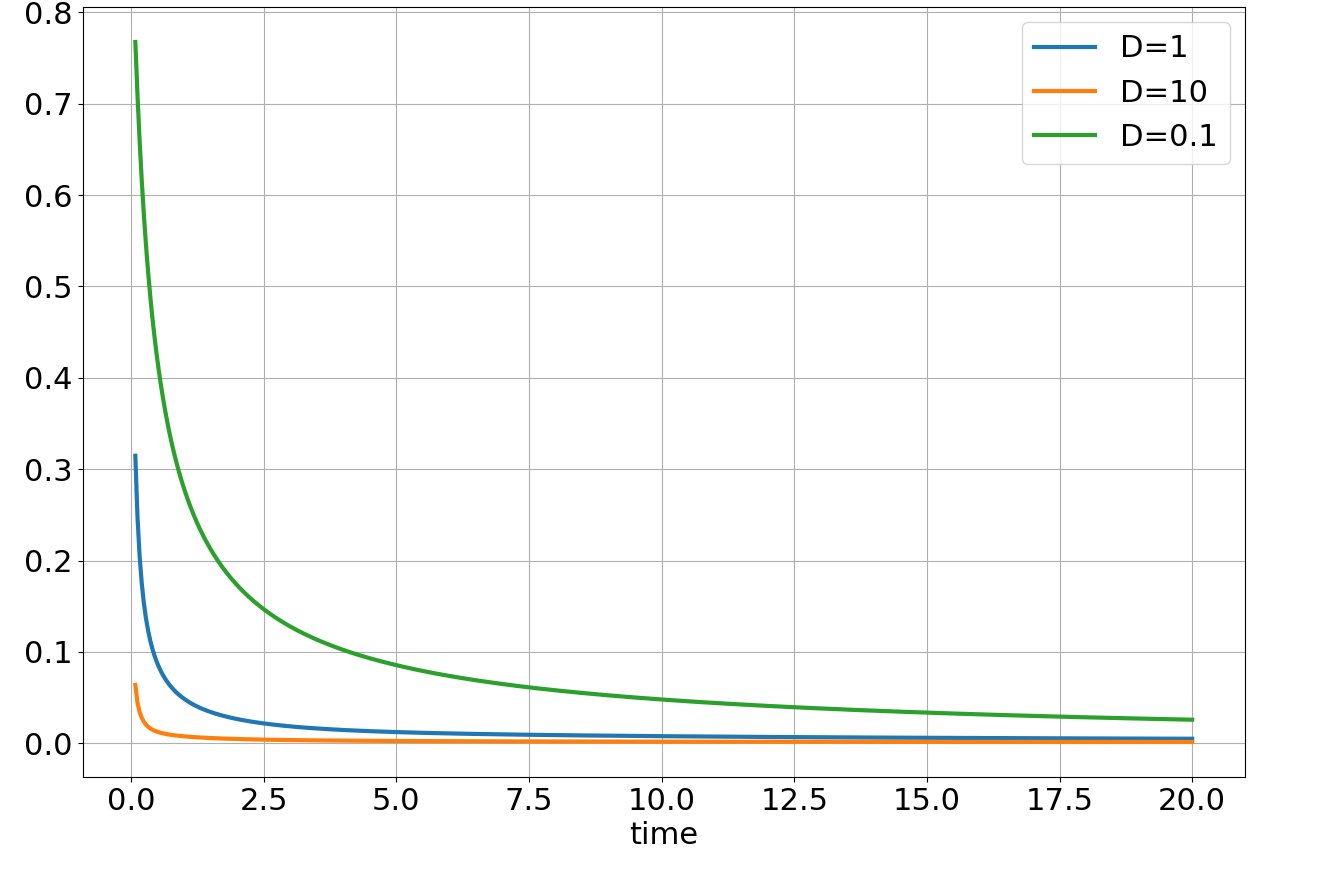}}
		\subfigure[Secretion and uptake: $\phi = a = 1.0$]{
			\includegraphics[width=0.4\textwidth]{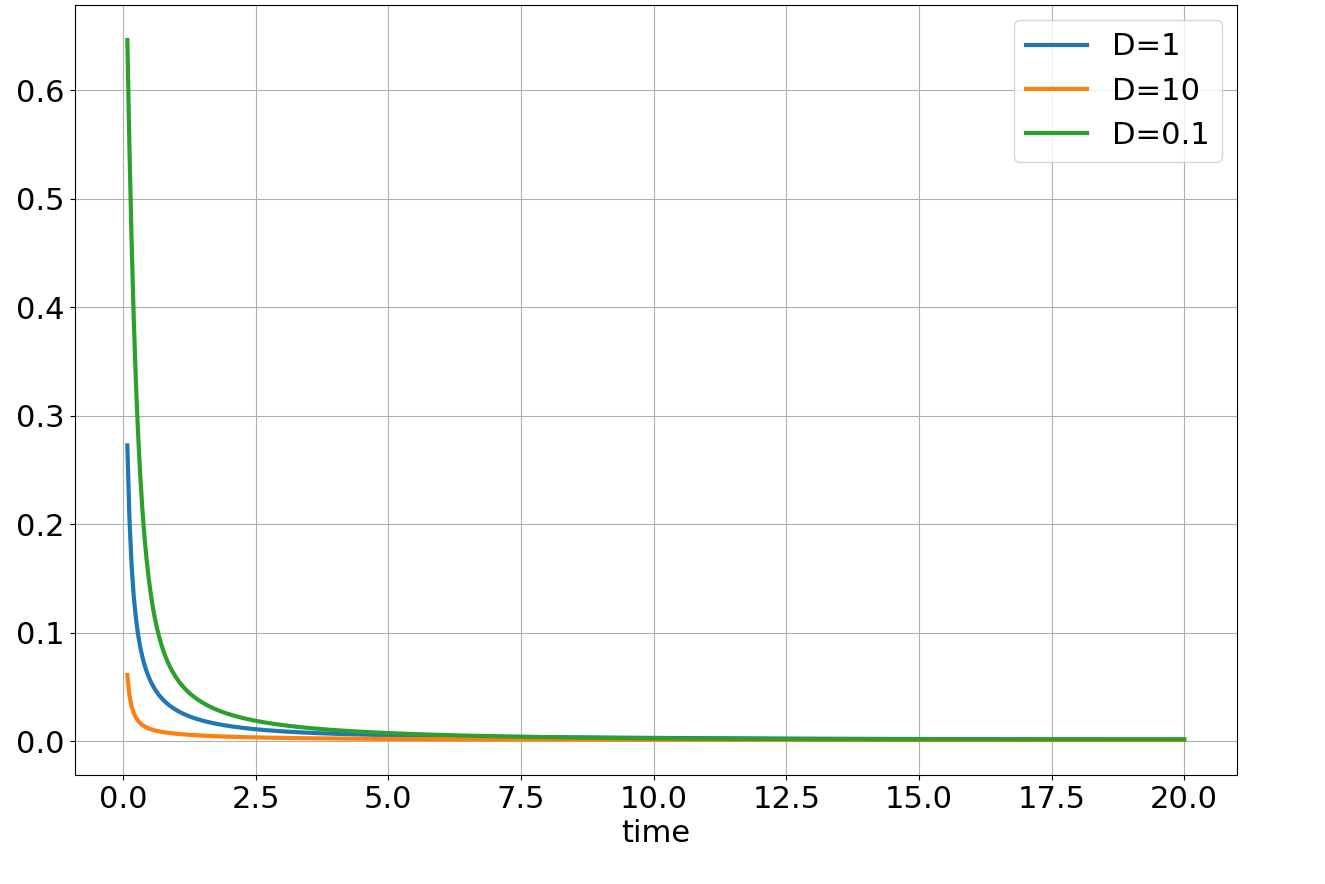}}
		\caption{Relative $L^2$-norm difference between the models, i.e. $e_{L^2(\tilde\Omega)}(t)$ given by \eqref{eq:rela_diff}, when the initial condition is zero in the computational domain. Here, diffusion coefficient $D$ is selected from $\{0.1, 1.0, 10.0\}$. Panel (a): cell only releases compounds to the ECM; (b) cell releases and takes up the compounds from the ECM. Different colors represent different values of $D$.}
		\label{fig:rela_diff_zero_IC}
	\end{figure}
	
	To investigate the influence of the difference between the solutions to the model after introducing the uptake of the compounds, we introduce a relative $L^2$-norm difference defined by
	\begin{equation}
	\label{eq:rela_diff}
	e_{L^2(\tilde\Omega)}(t) = \frac{\Vert u_P(\boldsymbol{x}, t) - u_S(\boldsymbol{x}, t)\Vert_{L^2(\tilde\Omega)}}{\Vert u_S(\boldsymbol{x},t)\Vert_{L^2(\tilde\Omega)}}.
	\end{equation}
	Furthermore, it is of interest to investigate the impact of different values of diffusion coefficient. Here, we consider $D \in \{0.1, 1.0, 10.0\}$. 
	Figure \ref{fig:rela_diff_zero_IC} shows the relative $L^2$-norm difference whether the uptake is taken into account or not. Aligning with the observations in \citet{Peng2023}, a larger diffusion coefficient will decrease the difference between the solutions to the two models, as it takes less time for the compounds to `travel' from the cell centre to the membrane in the point source model. Moreover, since there is no compound in the ECM initially, including the uptake does not cause a significant difference on $e_{L^2(\tilde\Omega)}(t)$ for small $t$. However, when the uptake is included in the model, $e_{L^2(\tilde\Omega)}$ decreases faster, even though in both cases, $e_{L^2(\tilde\Omega)}$ vanishes to almost zero. This is mainly because in the case of secretion and uptake, there is no compounds-exchanging anymore once $\phi - au = 0$, where $u$ is the concentration of the compounds. Note that this does not necessarily mean $u = 0$.  
	
	Comparison of the time scale of convergence to a steady level in Figure \ref{fig:steady} and Figure \ref{fig:rela_diff_zero_IC} suggests that the solutions to spatial exclusion model and point source model converge equally fast to steady state. Although the solution to the point source model differs initially from the other, on short time scale, it seems to capture equally well the convergence to the steady state on the longer time scale.
	
	
	
	\section{Conclusion}\label{sec:conclusion}
	Being interested in models for the exchange of compounds between (possibly moving) cells and the ECM, we modified the two models proposed in \cite{Peng2023}, leading to systems $\rm (BVP_S)$ and $\rm (BVP_P)$. The non-local character of the amplitude of the Dirac measure at the cell centre in the point source model makes analysis challenging. However, this choice of amplitude is strongly motivated by the case of emission only, that was studied in \citet{Peng2023} and that yielded good correspondence of solutions.
	
	We showed the well-posedness of the two models: the existence, uniqueness and continuous dependence on initial conditions of both models can be verified in suitable Sobolev spaces over the the respective computational domains. However, the Dirac measure used in the point source model causes the solutions to be in $W^{1,p}(\Omega)$ for $1<p<2$, but not in $H^1(\Omega)$ -- the space commonly used for FEM. It is expected to be in $H^1(\tilde\Omega)$, when similar argumentation as in \citet{HMEvers2015} is used, building on estimates for the heat kernel and its gradient on suitable bounded domains (with Neumann boundary conditions).

	The ultimate objective is to compare the solutions to the two models. Numerically, one can conclude that, after introducing the uptake, both system will reach a steady state in long term. In other words, the steady state appears when there is no exchange of compounds, which differ from the case when there is only secretion of the compounds. Furthermore, regardless of the value of the diffusion coefficient, the discrepancy between the solutions to the two models vanishes, shorter when there is exchange of compounds. 
	
	Questions that remain to investigate are, for example, how the decay rate of the relative difference depends on the uptake rate $a$, the cell radius or secretive flux density $\phi$. The derivation of theoretical estimates for the (absolute) $L^2$- or $H^1$-difference on $\tilde{\Omega}$ are much hindered so far by lack of applicable techniques that apply to obtaining estimations on a subdomain ($\tilde{\Omega}$) of the full computational domain ($\Omega$). Estimations for relative differences have not been obtained so far. Further analytic quantification of the discrepancy between the solutions remains a topic of forthcoming study, as it is an essential ingredient for rational application of the proposed approximation of spatial exclusion models by suitable point source models. 
	
	\color{black}
	
	
	
	\bibliography{MATHMOD_Conference/Latex_Full/ifacconf}

\end{document}